\documentclass{article}
\usepackage[utf8]{inputenc}

\title{Hamiltonian Stationary Lagrangian Surfaces with Non-Negative Gaussian Curvature in K\"ahler-Einstein Surfaces}
\author{Patrik Coulibaly \\
e-mail: cpatrik@math.ubc.ca}
\date{\today}

\usepackage{amsmath}
\usepackage{amssymb}
\usepackage{amsfonts}
\usepackage{graphicx}
\usepackage{enumitem}
\usepackage{musicography}
\usepackage[numbers]{natbib}
\usepackage{cite}
\usepackage{url}
\usepackage{amsthm}
\usepackage{commath}

\newtheorem{thm}{Theorem}[section]
\newtheorem*{thm*}{Theorem}

\newtheorem{lem}[thm]{Lemma}
\newtheorem{cor}[thm]{Corollary}

\theoremstyle{definition}

\theoremstyle{remark}
\newtheorem*{rem}{Remark}

\newtheoremstyle{bfnote}%
{}{}%
{}{}%
{\bfseries}{}%
{ }%
{\thmname{#1}\thmnumber{ #2}\thmnote{ (#3)}}

\theoremstyle{bfnote}
\newtheorem*{obs*}{Observation}

\providecommand{\R}{\mathbb{R}} 
\providecommand{\Z}{\mathbb{Z}} 
\providecommand{\C}{\mathbb{C}} 
\providecommand{\Proj}{\mathbb{P}} 
\providecommand{\Hyp}{\mathbb{H}} 

\providecommand{\rn}[1]{ \textup{\uppercase\expandafter{\romannumeral#1}}} 
\begin{document}

\maketitle

\begin{abstract}
   In this paper, we give some simple conditions under which a Hamiltonian stationary Lagrangian submanifold of a K\"ahler-Einstein manifold must have a Euclidean factor or be a fiber bundle over a circle. We also characterize the Hamiltonian stationary Lagrangian surfaces whose Gaussian curvature is non-negative and whose mean curvature vector is in some $L^p$ space when the ambient space is a simply connected complex space form.
\end{abstract}
\section{Introduction}\label{setion:Introduction}
Let $(M,g,J)$ be a K\"ahler manifold of complex dimension $n$. $M$ carries a natural symplectic structure given by the closed $2$-form $\omega$ which is defined by $\omega(X,Y) = g(JX,Y)$ for $X, Y \in T_pM$. We say that a Lagrangian submanifold $L \subset M$ is \textit{Hamiltonian stationary} if it is a critical point of the volume functional under compactly supported Hamiltonian deformations, i.e. variations for which the variational vector field is of the form $V=J\nabla f$ for some $f \in C^\infty_c(L)$. In \citep{oh}, Oh calculated the Euler-Lagrange equation of the variational problem and found that Hamiltonian stationary Lagrangian submanifolds are characterised by
\[
    \delta \alpha_H = 0,
\]
or equivalently by
\[
    \text{div}_{L}(JH) = 0,
\]
where $H$ denotes the mean curvature vector of $L$, which we define as the trace of its second fundamental form $A$, i.e. $H := Tr_g A$, $\alpha_H$ is the differential $1$-form on $L$ defined by $\alpha_H := \iota_{H} \omega = g(JH, \cdot)$ and $\delta$ is the co-differential operator on $L$ induced by the metric $g$.

By a theorem of Dazord (see, for example, Theorem 2.1 in \citep{oh}), in any K\"ahler manifold $M$, the restriction of the Ricci form $ric_{M}$ of $M$ to $L$ is given by $d \alpha_H$. When $M$ is K\"ahler-Einstein, i.e. $ric_{M} = c \omega$ for some constant $c$, then the differential $1$-form $\alpha_H$ is closed and thus defines a cohomology class in H$_{dR}^{1}(L) \cong H^1(L;\R)$ on any Lagrangian submanifold. Therefore, $\alpha_H$ is both closed and co-closed, hence harmonic,  on any Hamiltonian stationary Lagrangian submanifold of a K\"ahler-Einstein manifold.


In \citep{arsie2000maslov}, Arsie proved that when $M$ is a Calabi-Yau manifold, then  $\frac{1}{\pi} \alpha_H$ represents an integral cohomology class of $L$ called the \textit{Maslov class}. Therefore, we will refer to $\mu = \frac{1}{\pi} \alpha_H$ as the \textit{Maslov form} of $L$.

Let $NL$ denote the normal bundle of $L$ in $M$ and $\Gamma (NL)$ denote the collection of smooth sections of $NL$. Also, for any point $x \in L$ and vector $\overline{X} \in T_xM$, let $\overline{X}^\perp$ denote the projection of $\overline{X}$ onto $N_xL$ and let $\overline{\nabla}$ denote the Levi-Civita connection on $M$. Then, there is a connection $\nabla^\perp$ in $NL$ that is given by
\[
    \nabla^\perp_X V := \left(\overline{\nabla}_X V\right)^\perp
\]
for any normal vector field $V \in \Gamma(NL)$ and tangent vector $X \in T_xL$. We say that a normal vector field $V \in \Gamma(NL)$ is \textit{parallel} if $\nabla^\perp V\equiv  0$.

For any point $x \in L$ and vector $\overline{X} \in T_xM$, let $\overline{X}^\top$ denote the projection of $\overline{X}$ onto $T_xL$. Then, since $L$ is Lagrangian and $\nabla J = 0$, we have that
\begin{align*}
    \nabla_X JV = \left(\overline{\nabla}_X JV\right)^\top = \left(J\overline{ \nabla}_X V\right)^\top = J \left(\overline{ \nabla}_X V\right)^\perp = J\nabla^\perp_X V,
\end{align*}
for any normal vector field $V \in \Gamma(NL)$ and tangent vector $X \in T_xL$. Therefore, $JV$ is parallel if and only if $V$ is parallel. 

We say that $L$ has parallel second fundamental form if
\begin{align*}
    (\nabla_X A)(Y,Z) = \nabla^\perp_X A(Y,Z) -A(\nabla_X Y,Z) - A(Y, \nabla_X Z) 
\end{align*}
vanishes for all $X,Y,Z \in \Gamma(TL)$.

We present our results in two separate sections. In Section \ref{section:General}, we consider a complete, connected Hamiltonian stationary Lagrangian submanifold $L$ of arbitrary dimension $n$ inside a K\"ahler-Einstein manifold. We introduce a set of conditions, most of which consist of the non-negativity of the Ricci curvature of $L$ in the direction of $JH$ and some pointwise or integral control over the absolute value of $H$, that allows us to combine the Bochner formula for the harmonic $1$-form $\alpha_H$ and some Liouville-type theorems to deduce that $H$ must be parallel in the normal bundle of $L$. The existence of a non-trivial global parallel vector field can restrict both the topology and the geometry of a manifold significantly. For example, if $L$ is simply connected, then it must be isometric to a Riemannian product of the form $N \times \R$. As for a purely topological consequence, if $L$ is not diffeomorphic to such a product, then it must admit a circle action whose orbits are not homologous to zero. In Section \ref{section:surfaces_general}, we restrict our attention to the case when $n=2$. We also strengthen our assumptions by requiring that our surface has non-negative Gaussian curvature which allows us to describe explicitly all complete, connected Hamiltonian stationary Lagrangian surfaces in $\C^2$, $\C\Proj^2$ and in $\C\Hyp^2$ that has non-negative Gaussian curvature and whose mean curvature vector is in some $L^p$ space.

The author would like to express his gratitude to Prof. Jingyi Chen for his invaluable suggestions and support, which were essential to the completion of this paper.

\section{Hamiltonian Stationary Lagrangians in K\"ahler-Einstein Manifolds and the Bochner Method}\label{section:General}

Let $\mathcal{A}$ denote any of the following sets of assumptions:
\begin{enumerate}
    \item $Ric_L(JH,JH) \geq 0$ and $\abs{H} \in L^p$ for some  $p \in (2, \infty)$;
    \item $L$ has non-negative Ricci curvature and $\abs{H} \in L^p$ for some  $p \in (0, \infty)$;
    \item $L$ is oriented, $Ric_L(JH,JH) \geq 0$ and  $\abs{H} \rightarrow c := \inf_{L} \abs{H}$ as $r(x) \rightarrow \infty $  where $r(x) := d(x_0,x)$ is the distance function on $L$ relative to a fixed point $x_0 \in L$;
    \item $Ric_L(JH,JH) \geq 0$ and
    there exists a point $x_0 \in L$, a non-decreasing function $f:[0, \infty ) \mapsto [0, \infty )$, constants $C, R >0$ and $p \in (2,\infty)$ such that $\abs{H(x)} \leq f(r(x))$ for all $x \in L$ and 
    \begin{align}\label{ineq:log_growth}
    \frac{f(r)^pVol(B_r(x_0))}{r^2 \log(r)} \leq  C
    \end{align}
     whenever $r\geq R$. Here, $B_r(x_0)$ denotes the geodesic ball in $L$ of radius $r$ around the point $x_0 \in L$;
\item $L$ has conformal Maslov form, i.e. the vector field $JH$ is conformal.
\end{enumerate}

\noindent If at least one of the sets of assumptions labelled (1)--(5) is satisfied, we say that $\mathcal{A}$ is satisfied.

We can state the main result of this section as follows.

\begin{thm}\label{theorem:general}
     Let $L$ be a complete, connected Hamiltonian stationary Lagrangian submanifold of a K\"ahler-Einstein manifold. If $L$ satisfies $\mathcal{A}$, then
     \begin{itemize}
        \item[(a)] $H$ is parallel and thus has constant length;
        \item[(b)] $Ric_L(JH, \cdot)$ vanishes identically, so if there exists a point $x\in L$ such that $Ric_L|_x$ is non-degenerate, then $L$ must be minimal;
        \item[(c)] and the scalar curvature of $L$ must be constant along the integral curves of $JH$.
    \end{itemize}
\end{thm}

 The growth bound (\ref{ineq:log_growth}) from condition (4) is satisfied, for example, when $L$ has quadratic volume growth and $\abs{H}$ does not grow faster than $\log(r)^{\frac{1}{2p}}$ at infinity for some $1 < p < \infty$. In particular, it is satisfied when $L$ has quadratic volume growth and $\abs{H} \in L^\infty$. Therefore, we have the following corollary.

 \begin{cor}
       Let $L$ be a complete, connected Hamiltonian stationary Lagrangian submanifold of a K\"ahler-Einstein manifold. If $Ric(JH,JH) \geq 0$, $\abs{H} \in L^\infty$ and $L$ has quadratic volume growth, then the conclusions of Theorem \ref{theorem:general} hold; in particular, $\abs{H}$ must be constant.
 \end{cor}

\begin{rem}
      This phenomenon is related to the notion of parabolicity of a manifold. We say that a manifold is \textit{(strongly) parabolic} if it does not admit a negative, non-constant subharmonic function, i.e. if $f<0$ and $\Delta f \geq 0$, then it must be constant. It is easy to see that a parabolic manifold does not admit a non-constant subharmonic function that is bounded from above.  A sufficient condition\footnote{As it is discussed in \citep{grigor1999analytic} after Corollary 7.7, when the Ricci curvature is non-negative, then this condition is also necessary.} for the parabolicity of a manifold was given by Karp in \citep{karp1982subharmonic}, which implies, for example,  that every complete, non-compact manifold with quadratic volume growth is parabolic. 
\end{rem}

The main restriction imposed on $L$ by the conclusion of the Theorem \ref{theorem:general} is that $JH$ is parallel since the existence of a non-trivial global parallel vector field restricts the topology of a manifold significantly. For example, the following result of Welsh \citep{welsh1986existence}, states that the existence of a complete non-trivial global parallel vector field forces the existence of a circle action whose orbits are not real homologous to zero. By a complete vector field, we mean a vector field whose integral curves are defined for all time.

\begin{thm}[Welsh \citep{welsh1986existence}]\label{theorem:Welsh_non_cpct}
    Suppose that $M$ is a Riemannian manifold that admits a non-zero complete parallel vector field. Then either $M$ is diffeomorphic to the product of a Euclidean space with some other manifold, or else there is a circle action on $M$ whose orbits are not real homologous to zero. Moreover, if $M$ is not diffeomorphic to the product of a Euclidean space with some other manifold and its first integral homology class is finitely generated, then $M$ is a fiber bundle over a circle with finite structure group.
\end{thm}

Combining Theorem \ref{theorem:general} and Theorem \ref{theorem:Welsh_non_cpct} gives us the following corollary.

\begin{cor}\label{cor:circle_action}
    Let $L$ be a complete, connected Hamiltonian stationary Lagrangian submanifold of a K\"ahler-Einstein manifold. If $L$ is not minimal and it satisfies $\mathcal{A}$, then $L$ is diffeomorphic to the product of a Euclidean space with some other manifold or there is a circle action on $M$ whose orbits are not real homologous to zero. Moreover, it satisfies the conclusion of Theorem \ref{theorem:general}; and if it is not diffeomorphic to the product of a Euclidean space with some other manifold and its first integral homology class is finitely generated, then $M$ is a fiber bundle over a circle with finite structure group.
\end{cor}
    
\begin{proof}
    Suppose that $L$ is a complete, connected Hamiltonian stationary Lagrangian submanifold of a K\"ahler-Einstein manifold that satisfies $\mathcal{A}$. Then all the assumptions of Theorem \ref{theorem:general} are satisfied thus all of its conclusions hold. In particular, $JH$ is parallel so its integral curves are geodesics. Since $L$ is complete, all of its geodesics are defined for all $t \in \R$ and we see that $JH$ is a  complete vector field. Therefore, when $L$ is not minimal, $JH$ is a non-zero complete parallel vector field on $L$ and we can apply Theorem \ref{theorem:Welsh_non_cpct} to finish the proof.
\end{proof}

When $L$ is not minimal, we can use Corollary \ref{cor:circle_action} to establish the existence of a circle action on $L$ whose orbits are not real homologous to zero but only if $L$ is not diffeomorphic to a product of $\R$ and some other manifold. It turns out that when $L$ is simply connected then the existence of such a splitting is guaranteed. Moreover, $L$ can be split in such a  way isometrically.
\begin{cor}
    Let $L$ be a complete, connected Hamiltonian stationary Lagrangian submanifold of a K\"ahler-Einstein manifold. If $L$ is not minimal and it satisfies $\mathcal{A}$, then its universal cover $\pi: \Tilde{L} \to L$ equipped with the pull-back metric is isometric to $N \times \R$ for some totally geodesic submanifold $N$ of $\tilde{L}$.
\end{cor}
 
\begin{proof}
    Let $L$ be a complete, connected Hamiltonian stationary Lagrangian submanifold of a K\"ahler-Einstein manifold that satisfies $\mathcal{A}$. Suppose that $L$ is not minimal. Then, by Theorem \ref{theorem:general}, $JH$ is a non-zero parallel vector field.

Let  $\pi: \Tilde{L} \to L$ denote the universal cover of $L$ which we equip with the pull-back metric. This makes $\pi$ into a Riemannian covering. We know that $d \alpha_H = \delta \alpha_H = 0$ and since $\pi$ is a local isometry, we must also have $d \Tilde{\alpha}_H = \delta \Tilde{\alpha}_H = 0$ for $\Tilde{\alpha}_H = \pi^* \alpha_H$. Define $\widetilde{JH} = (\Tilde{\alpha}_H)^\sharp$. Then $JH = \pi_* \widetilde{JH}$ so $\widetilde{JH}$ must also be a parallel vector field. Since $\Tilde{L}$ is simply connected and $\Tilde{\alpha}_H$ is closed, there exists a smooth function $f \in C^\infty (\Tilde{L})$ such that $\Tilde{\alpha}_H = df$ or equivalently $\widetilde{JH} = \nabla f$. 

Since $\widetilde{JH} =\nabla f$ is parallel, by Lemma 2.3. in \citep{sakai1996riemannian}, $f$ is an affine function in the sense that  $f \circ \gamma: \R \to \R$ satisfies
\[
    f \circ \gamma (\lambda t_1 +(1-\lambda)t_2) = \lambda f \circ \gamma(t_1) + (1-\lambda) f \circ \gamma(t_2)
\]
for all maximal unit speed geodesics $\gamma$ in $\widetilde{L}$, $\lambda \in (0,1)$ and $t_1,t_2 \in \R$. Also, since $\widetilde{JH} =\nabla f$ is non-zero, $f$ is a non-trivial affine function so, by a theorem of Innami \citep{innami1982splitting}, $f^{-1}(0)$ is a totally geodesic submanifold of $\Tilde{L}$ and $f^{-1}(0) \times \R$ is isometric to $\Tilde{L}$. 
\end{proof}

In order to prove Theorem \ref{theorem:general}, we need the following lemma.

\begin{lem}\label{lem:step_1}
     Let $L$ be a complete, connected Hamiltonian stationary Lagrangian submanifold of a K\"ahler-Einstein manifold. If $L$ satisfies $\mathcal{A}$, then $\abs{H}$ is constant. 
\end{lem}
\begin{proof}
       Let $L$ be a complete, connected Hamiltonian stationary Lagrangian submanifold of a K\"ahler-Einstein manifold. Then $\alpha_H = (JH)^\flat$ is closed so we can apply the Bochner formula  \citep[p.~207]{petersen2006riemannian} for $\abs{\alpha_H}^2 = \abs{JH}^2 =\abs{H}^2$  to get
    \begin{align*}\label{eqn:bochner}
        \frac{1}{2} \Delta_L \abs{H}^2 &= \langle JH, \nabla \text{div}_{L} JH \rangle + Ric_L(JH,JH) + \abs{\nabla JH}^2  \\     
        &= Ric_L(JH,JH) + \abs{\nabla^\perp H}^2 
    \end{align*}

First, we observe that if $\mathcal{A}$ is one of the conditions (1)--(4)
, then $Ric_L(JH,JH) \geq 0$ and the function $\abs{H}^2$ is clearly subharmonic. Since compact manifolds do not admit non-constant subharmonic functions, we may assume without the loss of generality that $L$ is non-compact when  $\mathcal{A}$ is one of the conditions (1)--(4).

If $\abs{H}^2$ is in $L^p$ for some $1<p< \infty$, then by a well-known result of Yau \citep{yau1976some}, $\abs{H}^2$ is constant. Therefore, condition (1) implies that $\abs{H}$ is constant.

If $L$ has non-negative Ricci curvature, then by a result of Li and Schoen (Theorem 2.2. in \citep{li1984p}) it does not admit any non-negative $L^p$ subharmonic function for all $0 < p < \infty$. Thus, condition (2) implies that $\abs{H}$ is constant.

Now, assume that condition (3) is satisfied. In \citep{alias2019maximum}, Alías, Caminha and do Nascimento prove that every non-negative subharmonic function
that converges to 0 at infinity on a connected, oriented, complete and non-compact Riemannian manifold must be identically zero. Applying this maximum principle to the function $f = \abs{H}^2 -c$ gives us that $\abs{H}^2 \equiv c$. Therefore, we conclude that condition (3) also implies that $\abs{H}$ is constant.

 Next, assume that condition (4) is satisfied. Since $f$ is a non-negative and non-decreasing function,
\begin{align*}
    \frac{1}{r^2 \log r} \int_{B_r(x_0)} \abs{H}^p dV & \leq \frac{1}{r^2 \log r}  \int_{B_r(x_0)} f(r)^{p} dV\\
    & =\frac{f(r)^{p}Vol(B_r(x_0))}{r^2 \log r}  \\
    &\leq C
\end{align*}
whenever $r \geq R$. Therefore,
\[
    \limsup_{r\to \infty}{\frac{1}{r^2 \log r} \int_{B_r(x_0)} \left(\abs{H}^{2}\right)^{\frac{p}{2}} dV} < \infty
\]
However, in \citep{karp1982subharmonic}, Karp showed that every non-negative non-constant subharmonic function $g$ on a complete non-compact Riemannian manifold satisfies
     \[
        \limsup_{r\to \infty}{\frac{1}{r^2 \log r} \int_{B_r(x)} g^q dV} = \infty
     \]
for all $q \in (1, \infty)$ and center $x$. Therefore, $\abs{H}^{2}$ must be constant and we can conclude that condition (4) also implies that $\abs{H}$ is constant.

Finally, suppose that $JH$ is conformal. Then, since $JH$ is divergence-free,
\[
    \mathcal{L}_{JH} \hspace{2pt} g = \frac{2}{n} \text{div}(JH)g = 0.
\]
Therefore, the vector field $JH$ is in fact Killing and the tensor $\langle \nabla JH, \boldsymbol{\cdot} \hspace{2pt}\rangle$ is skew-symmetric. Since the dual $1$-form  $\alpha_H$ is closed, we also know that $\langle \nabla JH, \boldsymbol{\cdot} \hspace{2pt}\rangle$ is symmetric, and hence it must be zero. Therefore, $JH$ is parallel which implies that it must also have constant length. We can conclude that if $L$ has conformal Maslov class, then $\abs{H}$ must be constant, which completes the proof.
\end{proof}

 Now, we can prove Theorem \ref{theorem:general}.
\begin{proof}[Proof of Theorem \ref{theorem:general}] 
Let $L$ be a complete, connected Hamiltonian stationary Lagrangian submanifold of a K\"ahler-Einstein manifold that satisfies  $\mathcal{A}$. Then, as in the proof of Lemma \ref{lem:step_1}, we have that
\begin{align}\label{eqn:bochner}
     \frac{1}{2} \Delta_L \abs{H}^2 = Ric_L(JH,JH) + \abs{\nabla^\perp H}^2.
\end{align}
By Lemma \ref{lem:step_1}, $\abs{H}$ is constant so the left-hand side of equation (\ref{eqn:bochner}) vanishes identically. When $\mathcal{A}$ is one of the conditions (1)--(4), then we have two non-negative terms on the right-hand side so they must each vanish identically, i.e. we must have that $Ric_L(JH,JH) \equiv 0$ and $\abs{\nabla^\perp H}^2 \equiv 0$. When $JH$ is conformal, then by the same argument that we used in the proof of Lemma \ref{lem:step_1}, $JH$ is parallel. So $\abs{\nabla^\perp H}^2 \equiv \abs{\nabla JH}^2 \equiv 0$ which forces $Ric_L(JH,JH) \equiv 0$. Therefore, we can conclude that $H$ is parallel and $Ric_L(JH,JH)$ is identically zero whenever $\mathcal{A}$ is satisfied. This proves $(a)$.

Recalling the Weitzenb\"ock formula \citep[p.~211]{petersen2006riemannian}, we have
\begin{align*}
    \Delta \alpha_H = -Tr_g (\nabla^2 \alpha_H) + Ric_L(JH, \cdot)
\end{align*}
where $\Delta$ is the Hodge-Laplacian acting on differential $1$-forms. Since $\alpha_H$ is both harmonic and parallel, we have that  
\[
    \Delta \alpha_H = Tr_g (\nabla^2 \alpha_H) = 0
\]
and thus $Ric_L(JH, \cdot)$ must also vanish identically. Let us also assume that there exists a point $x \in L$ such that $Ric_L|_x$ is non-degenerate. Since  $Ric_L|_x(JH|_x,JH|_x) = 0$, we must have that $JH|_x = 0$. However, we know that $JH$ has constant length so $JH$ must vanish identically and thus $L$ is minimal. This proves $(b)$. 

Let us also recall the contracted Bianchi identity (Proposition 7.18. \citep{lee})
 \begin{align}\label{eqn:contracted_bianchi}
     \frac{1}{2}dS_L = Tr_g \nabla Ric_L
 \end{align}
where $S_L$ is the scalar curvature of $L$. The trace is taken on the first and the third indices, i.e. given a local orthonormal frame  $E_1, \dots, E_n$, equation (\ref{eqn:contracted_bianchi}) reads as

\begin{align*}
     \frac{1}{2} dS_L &= (\nabla_{E_i} Ric_L)( \cdot, E_i).
\end{align*}
Therefore, plugging $JH$ into equation (\ref{eqn:contracted_bianchi}) gives us that 
\begin{align*}
    \frac{1}{2} dS_L (JH) &= (\nabla_{E_i} Ric_L)(JH, E_i) \\
            &=E_iRic_L(JH,E_i) - Ric_L(\nabla_{E_i}JH, E_i) - Ric_L(JH, \nabla_{E_i}E_i) \\
            &=0.
\end{align*}
The first and the third terms vanish since $Ric_L(JH, \cdot) \equiv 0$, while the second term is zero because $JH$ is parallel. So we conclude that the scalar curvature must be constant along the integral curves of $JH$, which completes the proof of $(c)$.
\end{proof}

\section{Hamiltonian Stationary Lagrangian Surfaces with Non-Negative Gaussian Curvature in K\"ahler-Einstein surfaces}\label{section:surfaces_general}

Let $\Sigma$ denote a complete connected Hamiltonian stationary Lagrangian surface isometrically immersed in a K\"ahler-Einstein surface $M$.

In order to obtain a characterization that is more explicit than the one given by Theorem \ref{theorem:general}, we adjust our sets of assumptions from before. Let $\mathcal{A}'$ denote any of the following sets of assumptions:
\begin{enumerate}
    \item $\Sigma$ has non-negative Gaussian curvature and $\abs{H} \in L^p$ for some  $p \in (0, \infty)$;
    \item $\Sigma$ is oriented, it has non-negative Gaussian curvature and  $\abs{H} \rightarrow c := \inf_{\Sigma} \abs{H}^2$ as $r(x) \rightarrow \infty $;
    \item $\Sigma$ has non-negative Gaussian curvature and the growth condition (\ref{ineq:log_growth}) is satisfied;
    \item $\Sigma$ has non-negative Gaussian curvature and conformal Maslov form.
\end{enumerate}

\noindent If at least one of the sets of assumptions labelled (1)--(4) is satisfied, we say that $\mathcal{A}'$ is satisfied.

We will treat the cases when $\Sigma$ is compact and when it is non-compact separately.

\begin{thm}\label{theorem:kahler-einstein_cpct}
Let $M$ be a K\"ahler-Einstein surface and let $\Sigma$ be a closed connected Hamiltonian stationary Lagrangian surface in $M$ that satisfies $\mathcal{A}'$. If $\Sigma$ is orientable, then it is 
    \begin{itemize}
        \item a flat torus or
        \item a minimal sphere.
    \end{itemize}
If $\Sigma$ is not orientable, then it is 
\begin{itemize}
        \item  a flat Klein bottle or
        \item a minimal projective plane. 
    \end{itemize}
In both cases, $\Sigma$ has parallel mean curvature.
\end{thm}

\begin{thm}\label{theorem:kahler-einstein_non-cpct}
Let $M$ be a K\"ahler-Einstein surface. If $\Sigma$ is a complete, connected non-compact Hamiltonian stationary Lagrangian surface in $M$ satisfying $\mathcal{A}'$, then it has parallel mean curvature and it is
    \begin{itemize}
        \item isometric to $\R^2$,
        \item diffeomorphic to $\R^2$ and is minimal or
        \item it is flat and its fundamental group is isomorphic to $\Z$.   
    \end{itemize}

\end{thm}

\noindent We start by proving the following lemma.

\begin{lem}\label{lem:flat_or_minimal}
    Let $\Sigma$ be a complete, connected Hamiltonian stationary Lagrangian submanifold of a K\"ahler-Einstein surface $M$. If $\Sigma$ satisfies $\mathcal{A}'$ then it has parallel mean curvature and it is also flat or minimal.
\end{lem}

\begin{proof}
     Let $\Sigma$ be a complete, connected Hamiltonian stationary Lagrangian submanifold of a K\"ahler-Einstein surface $M$. Since dim$_\R (\Sigma)=2$, its curvature is entirely determined by its Gaussian curvature $K$ and, in particular, 
\[
    Ric_\Sigma = K g_{\Sigma}.
\]
Suppose that $\Sigma$ satisfies $\mathcal{A}'$. It is easy to see that $\mathcal{A}'$ is stronger than $\mathcal{A}$ so we can apply Theorem \ref{theorem:general} which tells us that $H$ is parallel and that
\begin{align*}
   K \abs{H}^2 \equiv  K \abs{JH}^2 \equiv Ric_\Sigma(JH,JH)\equiv 0.
\end{align*}
Since $H$ is parallel, it has constant norm and therefore $\Sigma$ must be minimal or flat.
\end{proof}
\begin{proof}[Proof of Theorem \ref{theorem:kahler-einstein_cpct}] 
     Let $\Sigma$ be a closed, connected Hamiltonian stationary Lagrangian submanifold of a K\"ahler-Einstein surface $M$. Also assume that $\Sigma$ satisfies $\mathcal{A}'$. Then, by Lemma \ref{lem:flat_or_minimal}, $\Sigma$ has parallel mean curvature and it must also be minimal or flat.

First, Suppose that $\Sigma$ is orientable. Then, by the Gauss-Bonnet theorem,
 \begin{align}\label{eqn:gauss-bonnet_cpct}
    \frac{1}{2 \pi} \int_{\Sigma} K dA = \chi(\Sigma)
 \end{align}
    where $\chi(\Sigma)$ is the Euler characteristic of $\Sigma$. Since $ \chi(\Sigma) = 2 - 
    2g$, where $g$ is the genus of $\Sigma$, and $K$ is non-negative, we see that the genus must be $0$ or $1$. Therefore, $\Sigma$ is diffeomorphic either to a sphere or to a torus respectively.  Equation (\ref{eqn:gauss-bonnet_cpct}) also tells us that $\Sigma$ is flat if and only if it has genus $1$, i.e. it is a torus. So, if $\Sigma$ is not flat then it is not just minimal but it must also have genus $0$ and thus it must be a minimal sphere. This completes the proof of the case when $\Sigma$ is orientable.

    Now, suppose that $\Sigma$ is not orientable. In this case, $ \chi(\Sigma) = 2 - \hat{g}$, where $\hat{g}$ is the non-orientable genus of $\Sigma$ which can be defined as the number of copies of $\R P^2$ appearing when the surface is represented as a connected sum of projective planes. Also, equation (\ref{eqn:gauss-bonnet_cpct}) still holds if we interpret the left-hand side as an integral of a density. One can easily see this by passing to the orientable double cover equipped with the pull-back metric. So, similarly to the orientable case, we have that $\hat{g}$ must be $1$ or $2$ and hence $\Sigma$ is diffeomorphic either to a real projective plane or to a Klein bottle respectively. We also see that $\Sigma$ is flat if and only if it is a Klein bottle. Therefore, when $\Sigma$ is not flat, then it is not just minimal but must also have non-orientable genus $1$ and thus it is a minimal real projective plane. This completes the proof of the non-orientable case. 
 \end{proof}

\begin{proof}[Proof of Theorem \ref{theorem:kahler-einstein_non-cpct}] 
It is known that the fundamental group of a non-compact surface is free (see, for example, \citep[p.~142]{stillwell2012classical}). The first singular homology group of $\Sigma$ with coefficients in $\Z$ is the abelianization of its fundamental group so $H_1(\Sigma;\Z)$ is the free abelian group on the generator set of $\pi_1(\Sigma)$. Since $H_1(\Sigma;\Z)$ is free, we know that  $H_1(\Sigma;\Z) = \Z^{b_1}$, where $b_1$ is the first Betti number of $\Sigma$. Therefore, the cardinality of the generator set of  $\pi_1(\Sigma)$ is equal to $b_1$.

First, assume that $\Sigma$ is orientable. Since $K \geq 0$, a result of Huber (Theorem 13. in \citep{huber1958subharmonic}) tells us that $\Sigma$ is finitely connected, i.e. it is homeomorphic to a closed surface with finitely many punctures. Therefore, $b_1$ must be finite. Moreover, since the top homology group of a non-compact manifold vanishes identically, we have that $b_2=0$ and thus $\chi(\Sigma) = 1 - b_1$. Also, by Theorem 10. in \citep{huber1958subharmonic},
\begin{align}\label{ineq:gauss-bonnet}
     \frac{1}{2 \pi} \int_\Sigma K  \leq  \chi(\Sigma)
\end{align}
so we have
\[
    b_1 \leq 1.
\]
If $\Sigma$ is not orientable, then applying the same argument but to the orientable double cover of $\Sigma$ also yields $ b_1 \leq 1$.

Since the generator set of $\pi_1(\Sigma)$ has either $0$ or $1$ element, $\pi_1(\Sigma)$ is either trivial or isomorphic to $\Z$. If $\pi_1(\Sigma) = \Z$, then $\chi(\Sigma) = 0$ so, by (\ref{ineq:gauss-bonnet}), $\Sigma$ must be flat.  If $\Sigma$ is simply connected, then it is diffeomorphic to $\R^2$. Finally, Lemma \ref{lem:flat_or_minimal} tells us that $\Sigma$ has parallel mean curvature and it is also minimal or flat which finishes the proof.
\end{proof}

\subsection{Hamiltonian Stationary Lagrangian Surfaces with Non-Negative Gaussian Curvature in Complex Space Forms}\label{subsection:surfaces_space_form}
Let $M(4c)$ be a complete, connected complex space form of complex dimension $2$ and constant holomorphic sectional curvature $4c$. Let $\Sigma \subset M(4c)$ be a Lagrangian submanifold.

We have the Wintgen-type inequality (Lemma 2.4. in \citep{MIHAI2014714}),
\begin{align}\label{ineq:Wingten_4c}
        K + \rho_N \leq c + \frac{\abs{H}^2}{4}
\end{align}

where $\rho_N \geq 0$ is a normalized (partial) normal scalar curvature\footnote{A similar inequality can be obtained using the Chen-Ricci inequality presented, for example, in \citep{deng2009improved}.}.

\noindent First, we look at the case $c=0$.
\begin{thm}\label{theorem:surface_0}
   Let $\Sigma \subset M(0)$ be a complete, connected Hamiltonian stationary Lagrangian submanifold. If $\mathcal{A}'$ is satisfied, then $\Sigma$ has parallel second fundamental form. Moreover, when the ambient manifold is $\C^2$, then $\Sigma$ is either
  \begin{itemize}
     \item a Lagrangian plane,
    \item a Riemannian product of a circle and a line (a Lagrangian cylinder), 
    \item or a Riemannian product of two circles (possibly of different radii).
\end{itemize}
\end{thm}
\begin{proof}
   Let $\Sigma \subset M(0)$ be a complete, connected Hamiltonian stationary Lagrangian submanifold. Then, by Lemma \ref{lem:flat_or_minimal}, $\Sigma$ has parallel mean curvature and it must also be minimal or flat. Also, by (\ref{ineq:Wingten_4c}),
\begin{align}\label{ineq:Wingten}
       0 \leq K \leq \frac{\abs{H}^2}{4}
\end{align}
so we see that if $\Sigma$ is minimal, then it must also be flat. Therefore, we may assume, without a loss of generality, that $\Sigma$ is flat. By Theorem 2.6. in \citep[p.~207]{kon1985structures}, $\abs{ \nabla A} = 0$ so we can conclude that $\Sigma$ has parallel second fundamental form.

For the rest of the proof, we assume that $M(0)$ is $\C^2$. Let $E_1, E_2$ be a local orthonormal frame on $\Sigma$. Then $E_1,E_2, JE_1, JE_2$ is a local orthonormal frame on $\C^2$ and the components $A_{ij}^k$ of the second fundamental form $A$ of $\Sigma$ in $\C^2$ are given by 
\[
    A(E_i,E_j)=A_{ij}^k JE_k.
\]
Let $A^k$ denote the $2 \times 2$ matrix $(A_{ij}^k)_{i,j}$. Then, since $\Sigma$ is flat, by Lemma 2.5. in \citep[p.~206]{kon1985structures}, its second fundamental form commutes, i.e. $A^kA^l = A^lA^k$ for all $k,l= 1,2$. 
Therefore,  by Theorem 2.9. in \citep[p.~210]{kon1985structures}, $\Sigma$ is congruent to one of the following standard Lagrangian submanifolds:
 \begin{enumerate}
    \item $\R^2 = \{(z_1,z_2) :\text{Im}(z_1)=0 \text{ and Im}(z_2)=0 \} \subset \C^2$ (a Lagrangian plane),
    \item $S^1(r) \times \R = \{(z_1,z_2) : \abs{z_1}^2 = r^2 \text{ and Im}(z_2)=0\}\subset \C^2$ for some $r>0$ (a Lagrangian cylinder),
    \item $S^1(r_1) \times S^1(r_2) = \{(z_1,z_2) : \abs{z_1}^2 = r_1^2 \text{ and }\abs{z_2}^2 = r_2^2\}\subset \C^2$ for some $r_1, r_2 > 0$ (a product of two circles).
   \end{enumerate}
\end{proof}

Before looking at the cases $c>0$ and $c<0$, we state some simple corollaries of Theorem \ref{theorem:surface_0}.
\begin{cor}
    Let $\Sigma \subset \C^2$ be a complete, connected Hamiltonian stationary Lagrangian submanifold that has non-negative Gaussian curvature and $\abs{H} \in L^p$ for some $0<p \leq \infty$. Then it is either a Lagrangian plane, a Lagrangian cylinder or a product of two circles. Moreover, it can only be a cylinder when $p = \infty$.
\end{cor}
\begin{proof}
    Let $\Sigma \subset \C^2$ be as stated in the corollary. When $p\in (0,\infty)$, then $\mathcal{A}'$ is clearly satisfied. Since $K \geq 0$, we know by the Bishop-Gromov volume comparison theorem that $Vol(B_r) \leq \pi r^2$ and thus $\Sigma$ has quadratic volume growth. Therefore, as discussed in the previous section after Theorem \ref{theorem:general}, the growth condition (\ref{ineq:log_growth}) is satisfied whenever $\abs{H} \in L^\infty$. So $\mathcal{A}'$ is satisfied when $p =\infty$ as well and we can use Theorem \ref{theorem:surface_0} to conclude that $\Sigma$ must be a Lagrangian plane, a Lagrangian cylinder or a product of two circles for any $0 < p \leq \infty$.

    Finally, we note that when $\Sigma$ is non-compact, then it has infinite volume \citep{yau1976some} so it must be minimal if it has a constant mean curvature that is in $L^p$ for some $p\in (0,\infty)$. The standard Lagrangian cylinder in $\C^2$ has constant mean curvature but it is neither compact nor minimal so it can only occur when $p=\infty$.
\end{proof}
We say that a complete non-compact submanifold $L$ is \textit{asymptoticaly minimal} if its mean curvature vector $H$ converges to $0$ at infinty, i.e. $\abs{H} \rightarrow 0$ as $r(x) \rightarrow \infty $.

\begin{cor}

    The only complete, connected, oriented and asymptotically minimal Hamiltonian stationary Lagrangian submanifolds of $\C^2$ with non-negative Gaussian curvature are Lagrangian planes.
\end{cor}

Since a complete K\"ahler manifold of positive holomorphic sectional curvature is necessarily simply connected (see, for example, \citep{tsukamoto1957kahlerian}), we may assume that $M(4)$ is $\C\Proj^2$ equipped with the standard Fubini-Study metric which has constant holomorphic sectional curvature $4$. Let $S^5 = \{ z \in \C^3 : \abs{z} = 1\}$ be the unit sphere in $\C^3$ equipped with induced metric. Then the map $\Pi: S^5 \to \C\Proj^2$ given by $x \mapsto [x]$, which is usually referred to as the Hopf fibration, can be used to construct Lagrangian immersions into $\C\Proj^2$. For more details, see, for example,  \S 3. in \citep{chen1997interaction}.
\begin{thm}
    Let $\Sigma \subset \C \Proj^2$ be a closed connected Hamiltonian stationary Lagrangian submanifold with non-negative Gaussian curvature. Then $\Sigma$ is
    \begin{itemize}
        \item a totally geodesic $\R \Proj^2$ or
        \item flat and is locally congruent to the image of $\Pi \circ L$ where $\Pi: S^5 \to \C\Proj^2$ is the Hopf fibration and $L: \Sigma \to S^5$ is given by $L(x,y) = (L_1(x,y), L_2(x,y), L_3(x,y))$ with
        \begin{align*}
            L_1(x,y) &= \frac{ae^{-i\frac{x}{a}}}{\sqrt{1 + a^2}}, \\ 
            L_2(x,y) &=\frac{e^{i(ax + by)}}{\sqrt{1 + a^2 + b^2}} \sin \left(\sqrt{1 + a^2 + b^2}y \right)\text{ and}\\
            L_3(x,y) &=\frac{e^{i(ax + by)}}{\sqrt{1 + a^2}} \left( \cos \left(\sqrt{1 + a^2 + b^2}y \right) - \frac{ib}{\sqrt{1 + a^2 + b^2}} \sin \left(\sqrt{1 + a^2 + b^2} y\right) \right)
        \end{align*}
          for some real constants $a \neq 0$ and $b$.
    \end{itemize}
\end{thm}
\begin{proof}
     Let $\Sigma \subset \C \Proj^2$ be a closed connected Hamiltonian stationary Lagrangian submanifold with non-negative Gaussian curvature. Then $\mathcal{A}'$ is satisfied so by Theorem \ref{theorem:kahler-einstein_cpct}, $\Sigma$ has parallel mean curvature and it is a minimal sphere, a minimal real projective plane, a flat Klein bottle or a flat torus.
     If $\Sigma$ is a minimal sphere or a minimal real projective plane then, by Theorem 7. in \citep{yau1974submanifolds}, $\Sigma$ is immersed in such a way that its image is a totally geodesic $\R\Proj^2$.
     If $\Sigma$ is a flat Klein bottle or a flat torus then $\Sigma$ has parallel second fundamental form by Theorem 2.6. in \citep[p.~207]{kon1985structures}. Therefore, the result follows from the classification of submanifolds with parallel second fundamental forms in $\C \Proj^2$  given by Theorem 7.1. in \citep{chen2010complete}.

\end{proof}

Finally, we consider the case $c < 0$. Let $\C\Hyp^2$ denote the complex hyperbolic space of constant holomorphic sectional curvature $-4$, let $\C^3_1$ denote $\C^3$ equipped with the psuedo-Euclidean metric $g = -dz_1d\Bar{z}_1 + dz_2d\Bar{z}_2 +dz_3d\Bar{z}_3$ and set $H^5_1 = \{ z \in \C^3 : g(z,z) = -1\}$.  Then the map $\Pi: H^5_1 \to \C\Hyp^2$ given by $x \mapsto [x]$, which we will also refer to as the Hopf fibration, can be used to construct Lagrangian immersions into $\C\Hyp^2$. For more details, see, for example,  \S 3. in \citep{chen1997interaction}.
\begin{thm}
     Let $\Sigma \subset M(4c)$ be a complete, connected Hamiltonian stationary Lagrangian submanifold for some $c < 0$. If $\mathcal{A}'$ is satisfied, then $\Sigma$ is flat and has parallel second fundamental form. Moreover, when the ambient manifold is $\C\Hyp^2$, then $\Sigma$ is locally congruent to the image of $\Pi \circ L$ where $\Pi: H^5_1 \to \C\Hyp^2$ is the Hopf fibration and $L: \Sigma \to H^5_1$ is given by $L(x,y) = (L_1(x,y), L_2(x,y), L_3(x,y))$ where
        \begin{enumerate}
            \item $L_1(x,y) =\frac{e^{i(ax + by)}}{\sqrt{1 - a^2}} \left( \cosh \left(\sqrt{1 - a^2 - b^2}y \right) - \frac{ib}{\sqrt{1 - a^2 - b^2}} \sinh \left(\sqrt{1 - a^2 - b^2} y\right) \right)$, \newline
            $L_2(x,y) =\frac{e^{i(ax + by)}}{\sqrt{1 - a^2 - b^2}} \sinh \left(\sqrt{1 - a^2 - b^2}y \right)$ and \newline
             $L_3(x,y) = \frac{ae^{i\frac{x}{a}}}{\sqrt{1 - a^2}}$ 
        for some real constants $a$ and $b$ satisfying $a \neq 0$ and $a^2 + b^2 < 1$;
         \item $ L_1(x,y) = \left( \frac{i}{b} + y \right)  e^{i(\sqrt{1-b^2} x + by)}$,\newline
            $L_2(x,y) = y e^{i(\sqrt{1-b^2} x + by)}$ and \newline
            $L_3(x,y) = \frac{ \sqrt{1 - b^2}}{b} e^{i\frac{x}{\sqrt{1-b^2}}}$ 
        for a real number $0 < b^2 < 1$;
          \item $L_1(x,y) =\frac{e^{i(ax + by)}}{\sqrt{1 - a^2}} \left( \cos \left(\sqrt{ a^2 + b^2 - 1}y \right) - \frac{ib}{\sqrt{a^2 + b^2 - 1}} \sin \left(\sqrt{a^2 + b^2 - 1} y\right) \right)$, \newline
            $L_2(x,y) =\frac{e^{i(ax + by)}}{\sqrt{a^2 + b^2 - 1}} \sin \left(\sqrt{a^2 + b^2 - 1}y \right)$ and \newline
             $L_3(x,y) = \frac{ae^{i\frac{x}{a}}}{\sqrt{1 - a^2}}$ 
        for some real constants $a$ and $b$ satisfying $0 < a^2 < 1$ and $a^2 + b^2 > 1$;
         \item $L_1(x,y) = \frac{ae^{i\frac{x}{a}}}{\sqrt{a^2 - 1}}$, \newline
            $L_2(x,y) =\frac{e^{i(ax + by)}}{\sqrt{a^2 + b^2 - 1}} \sin \left(\sqrt{a^2 + b^2 - 1}y \right)$ and \newline
           $L_3(x,y) =\frac{e^{i(ax + by)}}{\sqrt{a^2 - 1}} \left( \cos \left(\sqrt{ a^2 + b^2 - 1}y \right) - \frac{ib}{\sqrt{a^2 + b^2 - 1}} \sin \left(\sqrt{a^2 + b^2 - 1} y\right) \right)$   
        for some real constants $a$ and $b$ satisfying $a^2 > 1$;
        \item $L_1(x,y) = \frac{e^{ix}}{8b^2} \left( i + 8b^2(i + x) - 4by\right)$, \newline
        $L_2(x,y) = \frac{e^{ix}}{8b^2} \left( i + 8b^2x -4by\right)$ and \newline
         $L_3(x,y) = \frac{e^{i(x +2by)}}{2b}$
         for a real number $b \neq 0$; or 
          \item $L_1(x,y) = e^{ix}\left( 1 + \frac{y^
          2}{2} - ix\right)$, \newline
        $L_2(x,y) = e^{ix} y$ and \newline
         $L_3(x,y) = e^{ix}\left(\frac{y^
          2}{2} - ix\right)$.
        \end{enumerate}
\end{thm}
\begin{proof}
    Let $\Sigma \subset M(4c)$ be a complete, connected Hamiltonian stationary Lagrangian submanifold for some $c<0$. Suppose also that  $\mathcal{A}'$ is satisfied. Then by Lemma \ref{lem:flat_or_minimal}, $\Sigma$ has parallel mean curvature and it is also flat or minimal. However, since $K \geq 0$, it is clear from (\ref{ineq:Wingten_4c}) that $\Sigma$ cannot be minimal. Therefore, $\Sigma$ must be flat and thus it has parallel second fundamental form by Theorem 2.6. in \citep[p.~207]{kon1985structures}. When $M(4c)$ is $\C \Hyp^2$, the result follows from the classification of submanifolds with parallel second fundamental forms in $\C \Hyp^2$  given by Theorem 7.2. in \citep{chen2010complete}.
\end{proof}

\bibliography{biblio_arxiv}
\end{document}